\documentclass[a4paper,10pt]{article}
\usepackage[utf8]{inputenc}

\usepackage{graphicx}
\usepackage{amsmath}
\usepackage{amssymb}
\usepackage{amsthm}
\title{The replacement property of PSL$(2,p)$ and PSL$(2,p^2)$}
\author{Hy P.G Lam\footnote{This material is based on research projects supervised at the 2017 Cornell University Math SPUR program. The author of this paper is grateful for the academic support and guidance from Professor R. Keith Dennis and graduate mentor Ravi Fernando.}\\University of California, Berkeley\\Berkeley, CA 94720}
\date{September 20th, 2017}

\begin{document}

\maketitle

\begin{abstract}
    In 2014, Benjamin Nachman \cite{nachman} showed that when $p\equiv$1 mod 8,  the 2-dimensional projective linear group over the field of $p$ elements fails the replacement property if the maximal length $m$ of an irredundant generating sequence for the group is 3. In addition, if $m=4$, the group satisfies the property for any prime $p$. In this paper, we will extend such classification for PSL$(2,q)$ where $q$ is $p$ and $p^2$ with conditions of modulo 8 and 10 on $p$. 
\end{abstract}

\section{Introduction}
Given an arbitrary group $G$, we denote $s=(g_1,...,g_k)$ a finite sequence of elements $g_i$'s in $G$. $s$ is said to be an \emph{irredundant generating sequence} for $G$ if $g_1,...,g_k$ generate $G$ and, for any $i\in\{1,...,k\}$, $g_i \notin \langle g_j:j\neq i\rangle$ \footnote{$s$ is also said to be \emph{irredundant} or \emph{independent} if the second condition holds.}. We note that any finite generating sequence always contains an irredundant generating one since we can just remove generators that can be formed by the others until the process terminates. Also, if $G$ is finite, there is an irredundant generating sequence in $G$. However, there may not exist any finite irredundant generating sequence in the case of infinite group. An example is the additive group $\mathbb{Q}$ of rational numbers.\\ From here\footnote{For the interest of this paper, we restrict our attention to only finite groups.}, we see that the idea of an irredundant generating sequence is analogous to that of a basis for a vector space, where, of course, irredundance is the general quality of linear independence and the generating property refers to the span of the basis vectors to the entire space. Likewise, the \emph{replacement property} also arises naturally as a notion from linear algebra generalized for arbitrary finite groups. That is, for a finite dimensional vector space with a given basis, any nontrivial vector can replace a basis vector to form a new basis. In fact, for any nontrivial linearly independent set $I$, there exists a subset of same cardinality of the basis that can be replaced by $I$ to give a new basis. Such intuition allows R.Keith Dennis and Dan Collins to arrive at the following definition.\\ \textbf{Definition}: Given a length-$k$ irredundant generating sequence $s$=$(s_1,...,s_k)$ for a finite group $G$, $s$ is said to \emph{satisfy the replacement property} if for any nontrivial element $g\in G$, there is a slot $i$-th in $s$ so that $g$ can replace $g_i$ to give a new generating (not necessary independent) sequence for $G$. Furthermore, $G$ is said to \emph{satisfy the replacement property for length-$k$} if all irredundant generating sequences of length $k$ satisfy the replacement property. Lastly, if we replace $k$ with $m$ (formally noted as $m(G)$), where $m$ denotes the maximal length of an irredundant generating sequence, $G$ \emph{satisfies the replacement property} if $G$ satisfies the replacement property for length-$m$. 
We sometimes abbreviate the replacement property as RP for shorter notation. 

Since the definition originates from elementary properties of finite dimensional vector spaces, it is trivial that any such vector space satisfies the replacement property. However, the quaternion group $Q_8$ mentioned in Nachman's paper is a classic example where this property is not satisfied in general. In fact, the author shows there exists an infinite class of simple groups, namely PSL$(2,p)$ where $p\equiv$ 1 mod 8 and $m(\mathrm{PSL}(2,p)) =3$, that fails the property \cite{nachman}. The proof serves as a template for us to arrive at similar results, which provide a larger set of cases where the replacement property fails to hold in general for PSL$(2,q)$. 
Before we begin our discussions about these cases, the motivation for studying the replacement property of a given finite group comes from our goals to generalize the study of linear algebra to obtain deeper understandings of finite group structures. In the attempt to answer the question whether a given group satisfies RP or not, we gain a clearer intuition about its subgroup structures as well as the characteristics of the generators in terms of the subgroups they generate. However, it requires powerful computing equipments with large storage of data to fully describe all possible subgroup structures for an arbitrary finite group, let alone their interactions. Classic examples are the sporadic simple groups: Mathieu groups $M_{11},M_{12},M_{22},M_{23},M_{24}$, Janko groups $J_1,J_2,J_3,J_4$, Baby Monster group $F_2$ and the Monster group $F_1$ etc.\cite{wilson}, although there are several remarkable constructions that provide combinatorial structures to these gigantic groups (see \cite{graphtower}). Fortunately, being a class of simple groups that is fundamental to the study of classical groups, the 2-dimensional projective linear group over the field of $q$ elements enjoys a nice maximal-subgroup description due to the work of L.E Dickson in 1901\cite{dickson}. In addition, we recall Galois's construction of PSL$(2,q)$ from the general linear space GL$(2,q)$, where PGL$(2,q):=\mathrm{GL}(2,q)/\{\alpha \mathbb{I}:\alpha \in F^{*}_q \}$ is isomorphic to the group of all Mobius transformation from $P^1(F_q)=F_q \cup \{\infty\}$ to itself via the natural isomorphism $\rho: \begin{pmatrix} a &b \\c &d \end{pmatrix} \mapsto (z\mapsto \frac{az+b}{cz+d})$. The map $\rho|_{\mathrm{PSL}(2,q)}$ is then an isomorphism from PSL$(2,q)$ to the group of Mobius transformations whose determinant is a square in $F_q^{*}$. Such congruent structure gives a useful interpretation of the subgroups in PSL$(2,q)$ in terms of their actions on $P^1(F_q)$.\\\\\textbf{Theorem}\cite{king}: Given a maximal subgroup in PSL$(2,q)$, it is an isomorphic copy of the following classes: 
\begin{itemize}
\item A group of order $q(q-1)$ i.e $C_q\rtimes C_{(q-1)/2}$ stabilising a point on $P^1(F_q)$
\item Dihedral group $D_{q-1}$ for $q$ odd, $D_{2(q-1)}$ for $q$ even, both are pair stabilisers on $P^1(F_q)$.
\item Dihedral group $D_{q+1}$ for $q$ odd, $D_{2(q+1)}$ for $q$ even, both are pair stabilisers on $P^1(F_q)$.
\item Subfield group PSL$(2,q_1)$ for $q$ an odd prime power of $q_1$, PGL$(2,q_1)$ for $q=q_1^2$, $q$ odd 
\item $S_4$ if $q=p\equiv \pm$ 1 mod 8, or $q=p^2$ with $p\equiv \pm$ 3 mod 8 ($p>3$). 
\item $A_4$ if $q=p\equiv \pm$ 3 mod 8 ($p>3$). 
\item $A_5$ if $q=p\equiv \pm$ 1 mod 10, or $q = p^2$ with $p\equiv \pm$ 3 mod 10. 
\end{itemize}

\section{Families of $\mathrm{PSL}(2,q)$ that fail RP}\quad\\ \textbf{Theorem 2.1}:          
If $p\equiv \pm$ 3 mod 8 ($p$ prime, $p \geq 7$), $G$ = PSL$(2,p^2)$ fails the replacement property if $m(G)=3$. 

The condition for $p \equiv \pm$ 3 mod 8 is important for the proof of this statement since it allows us to utilize Dickson's theorem on the classification of maximal subgroups of PSL$(2,q)$ in which $S_4$ is a maximal subgroup if $q=p^2$. Furthermore, the schematic approach to the proof is analogous
to the one introduced by B. Nachman[3]. Specifically, one constructs an irredundant generating sequence  of length 3 and an element in $G$, namely the 90-degree rotation matrix. Such element once replaces any of the generating elements, fails to give a new generating sequence since the subgroup generated by the new sequence is contained in a maximal subgroup in PSL$(2,p^2)$.\\\\ \textbf{Proposition 2.2}:  Any quadratic
polynomial over prime field $F_p$ splits $F_{p^2}$ ($p$ prime). Thus, equations such as $x^2 -2$ always have
a solution in $F_{p^2}$. 

\begin{proof}
Suppose given $f(x)$ a polynomial over $F_p$ of degree 2, if f is reducible,
we are done. If not, $F_p[x]/{ \langle f(x)\rangle}$ is a field since irreducibility
implies
$\langle f(x)\rangle$ being prime in principal ideal domain, and
being prime
implies maximality in unique factorization domain. Also,  the field  has
$p^2$ elements because, for any element $g(x) $ in  $F_p[x]$, there exist
$q(x)$, $r(x)$ polynomials in  $F_p[x]$ such that $g(x) = q(x)f(x) + r(x),$
where $deg(r) \leq 1,$ and thus, there are such $p^2$ possibilities for $r(x)$.
It is a standard result that any finite splitting field of a given order
is unique up to isomorphism. Therefore, $f$ over $F_p$ splits  $F_{p^2}\cong F_p[x]/{\langle f(x) \rangle}$ and, thus, has a root in $F_{p^2}$.\end{proof}

We mention another simple, yet useful proposition.\\

\begin{noindent}\end{noindent}\textbf{Proposition 2.3}: Any traceless element in SL$(2,F)$ has order 4. If any element has
trace 1, it has order 6.

\begin{proof}It is straightforward by matrix multiplication and definition
of the group.\end{proof}
\begin{proof}(\emph{of theorem 2.1}). Consider the following matrices in
    $SL(2,p)$

 $M$ = 
$\left.\begin{pmatrix} x &y\\y & -x \end{pmatrix}\right.$
,$N$=
$\left.\begin{pmatrix} z &t\\t & -z \end{pmatrix}\right.$, 
$W$ = $\left.\begin{pmatrix} 0&-1\\1&0\end{pmatrix}\right.$, where the entries are elements in $F_{p^2}$.
By previous proposition, we have $M , N,W $ are of order 4. Thus, their images
under the canonical projection from $SL(2,p^2)$
to $PSL(2,p^2)$ have order 2. We let $m,n,w$ be those corresponding images.
By the same reason, we also have $mw , nw $ having order 2. Thus, $\langle m,w \rangle$ $\cong$ $\langle n,w\rangle$ $\cong$ $K_4$.
This means $M$,$N$ are not $W$. 

Now, we define $R=\begin{pmatrix} a_1 & 0\\-1 & -a_1\end{pmatrix}$ in $SL(2,p^2)$. Since
$-a_1^2 =1 $, $a_1 $ is of order 4 in the multiplicative field $F_{p^2}^*$.
Such an element exists because  $p^2-1$ is divisible by 8, hence divisible by 4. This shows that $R$ is well-defined. We have that $WR=\begin{pmatrix}
1 & a_1\\a_1 & 0\end{pmatrix}$ , thus $wr$ has order 3, where $r$ is the
projection of $R$ in $PSL(2,p^2)$. In addition, since $r$ and $w$ have order
2, by writing $rw=w^{-1}(wr)w,$ $1=wrrw$ and $r=(rw)w$,      we obtain $\langle r,w \rangle$=$\langle wr,w\rangle$=$\langle wr,w:w^{-1}(wr)w=(wr)^{-1},
(wr)^{3}=1=w^2\rangle$ $\cong$ $S_3$
noncyclic, which means $r,w$ are distinct.  
    
\textbf{Claim 2.3.1}: $\langle wm,wr\rangle\cong S_4$, so is $\langle
wn,wr\rangle$.
      
      \emph{Proof of claim.}
      Since $WM=\begin{pmatrix} -x & y\\y &x\end{pmatrix}, 
      WR=\begin{pmatrix}1 &a_1 \\a_1& 0 \end{pmatrix},
      WMWR=\begin{pmatrix}-x+a_1y&y+a_1x \\-a_1x & a_1y\end{pmatrix}$, by
Fricke's lemma [2],
\begin{align*}
&Tr([WM,WR])\\
&\quad=Tr(WM)^2+Tr(WR)^2+Tr(WMWR)^2\\
&\qquad-Tr(WM)Tr(WR)Tr(WMWR)-2\\
&\quad=2x^2-4a_1xy -3y^2
\end{align*} (using the fact that $x^2+y^2 =-1)$.
By D.McCullough [1], if we set this trace value of the commutator of $WM,
WR$ to be 1, then the projection of the subgroup generated by $WM,WR$ is
isomorphic to $S_4$. From this, it is left for us to show that there exist
elements $x,y$ in $F_{p^2}$ that satisfy such constraint. By solving the
quadratic equation $3x^2 -4a_1xy-2y^2=0$, we obtain $y=(-a_1+\frac{1}{\sqrt{2}})x
\Rightarrow  x^2=\frac{\sqrt{2}\pm2a_1)^2}{9}$. Since $p^2\geq 49$, $3^{-1}$
is well-defined in $F_{p^2}$. By lemma 1, $\pm \sqrt2$ is also well-defined.
Thus, such element $x,y$ exist in $F_{p^2}$, so do $z,t$. In fact, since
there are at least two distinct square root of 2 in $F_{p^2}$, the 2 tuples
$(x,y), (z,t) $ are uniquely determined, which implies $M,N $ are distinct, hence the claim.

\textbf{Claim 2.3.2}: $\langle wm,wr\rangle=\langle m,r,w\rangle$, likewise,
$\langle wn,wr\rangle=\langle n,r,w \rangle$\

\emph{Proof of claim.} We have: $w(mw)r=w(wm)r=mr, (mr)(rm)=1$, which implies
$mr,rm \in\langle wm,wr \rangle$. Given any even-length word comprising of
$m,r,w$, since
$r^{-1}=r, w^{-1}=w,m^{-1}=m                        ,wm=mw,(wr)^{-1}=rw $,
such word\footnote{Since $r,w,m$ are of order 2, such word is assumed to
be written with nonnegative powers.} is contained in $\langle wm,wr \rangle$.
Now,
for any odd-length word written by $m,r,w$, it is a word obtained from an
even-length word with either $m,r$ ,or ,$w$. Therefore, we have $\langle m,r,w\rangle\subset\bigcup_{j\in\{m,r,w,id\}}j.\langle wm,wr\rangle \Rightarrow|\langle m,r,w\rangle|\leq4|\langle wm,wr\rangle|=4|S_4| <|G|$
 (since $p^2>49$). This shows that $\langle m,r,w\rangle$ is a proper
subgroup of $G$ containing $\langle wm,wr\rangle$. By
claim  2.3.1, $\langle m,r,w\rangle=\langle wm,wr\rangle$. By symmetry, the desired
result is also true for
$\langle wn,wr\rangle$.

\textbf{Claim 2.3.3}: $\langle wm, wn\rangle$ is proper. 

\emph{Proof of claim}.  Since both $wm,wn$ are of order 2, we have $\langle
wm,wn\rangle
\cong D_{2k}$, where $k=ord(wmwn)=ord(mn)\geq 2 $, which is non-abelian.
Since $G$ contains at least 2 copies of $K_4$, thereby, not being dihedral.

\textbf{Claim 2.3.4} :  $\langle m,n,w\rangle $ is also proper. 

\emph{Proof of claim}.  Since one can write $w(mw)n=w(wm)n=mn \in\langle
wm,wn\rangle$
and $m,n,w$ have order 2, as in the proof of claim 2.3.2, we can apply the same
argument about the even and odd-length words written by $m,n,w$ to arrive
at the inequality $|\langle m,n,w\rangle| \leq4|\langle wm,wn\rangle|$. By
Dickson's classification
of  subgroups of $PSL(2,q)$ [4], $G$ contains the largest dihedral  subgroups
of order $q+1$. For any $p>7$, $4(p^2 +1)\leq 2^{-1}(p^2-1)p^2(p^2+1)\leq|G|$.
This and the previous inequality give the claim.

\textbf{Claim 2.3.5}:  $(wm,wn,wr)$ is an irredundant generating
sequence of
      $G$ and $w $ cannot replace any of the generating elements, which imply
the main result. 

\emph{Proof of claim}. We must have $ wn\notin\langle wm,wr\rangle$. Otherwise,
by claims
2.3.1 and 2.3.2, $\langle m,w,r\rangle=\langle n,w,r\rangle\Rightarrow mn\in S_4$. However,
we notice that\\

$Tr(MN) =2(xz+yt)= 
\begin{cases} \frac{1}{3}{(\sqrt{2} +a_1)}^2  & \text{if $x=(-a_1+\frac
{1}{2})y,y=\frac{1}{3}(2a_1-\sqrt{2})$} \\
\frac{1}{3}({\sqrt{2} -a_1})^2 & \text{if $x=(-a_1 -\frac{1}{2})y, y=\frac{1}{3}(2a_1+\sqrt{2})$}
\end{cases}$\\

Also, for $p^2>49$, it is straightforward to check that $\frac{2}{3}(\sqrt{2}
\pm a_1)^2 \neq \pm \sqrt{2}$. From this, it is clear to see $mn$ has order
strictly greater than 4, which is impossible to be contained in $S_4$. By
the same reasoning, $wm\notin \langle wn,wr\rangle$. If $wr\in \langle wn,wm\rangle$,
we have $\langle wm,wr\rangle
\cong S_4 \lneq\langle wn,wm \rangle $, which is impossible by claim 2.3.3 and
maximality of
$S_4$. This proves the irredundence of the sequence. 

By the previous paragraph, $\langle wm,wn,wr\rangle$ properly contains $\langle
wm,wr\rangle $, which
is maximal. We deduce that the sequence generates $G$. 

Now, if we consider $\langle w,wn,wr\rangle  $ or $\langle wm,w,wr\rangle$,
by claims 2.3.2 and 2.3.3, these subgroups
are proper. Since $\langle wm,wn,w \rangle $ is a subgroup of $\langle m,n,w\rangle$,
by claim 2.3.4, this
subgroup is proper. As a result, the sequence fails the replacement property.
\end{proof}

\begin{noindent}\end{noindent}\textbf{Theorem 2.4}: For $G$=PSL$(2,p)$, for prime $p=\pm1$ mod $10$ and $1$ mod $4$, $G$ fails the replacement property if $m=3$. \

For the proof of theorem 2.4, we follow a similar approach as done for theorem 2.1. That is, we construct a length-3 irredundant generating sequence $s$ for $G$ and an element in $G$ that cannot replace any generator in $s$ to give a new generating sequence. However, subtle details regarding maximal subgroups and orders of the constructed elements are to be addressed in order to arrive at the desired result. We first start with a proposition.\\

\begin{noindent}\end{noindent}\textbf{Proposition 2.5}: For $p>5$, any elements in SL$(2,p)$ possessing trace, whose value is a root of the equation $t^2 -t -1$ defined in $F_p$, has order 5. 
\begin{proof}
 We recall that the $m$-th cyclotomic polynomial $\Phi_m$ is the unique polynomial over an arbitrary field $F$ that divides only $x^m -1$ and is not a divisor of $x^r-1$ for any $r <m$, whose roots are all the $m$-th roots of unity, namely,
$\Phi_m(x) = \prod_{1\leq k\leq m:gcd(k,m)=1}(x-x^ {\frac{k}{m}})$. Also, in order for an element $X$ in SL$(2,F)$ to have have order $m$, the eigenvalues $\lambda$, $\lambda^{-1}$ of $X$ are roots of $\Phi_m(x)$ (assuming $m>2$). This happens iff the characteristic polynomial of degree 2 divides $x^m-1$. Now, suppose $X$ is an element in SL$(2,F_p)$ whose trace is a root of the equation $x^2 -x-1$, that is $Tr(X)= \frac{1\pm \sqrt{5}}{2}$. The characteristic polynomial $p(x)$ of $X$ is of form $x^2 + tx +1$ because the constant term is the determinant of $X$, which is 1 and $t= Tr(X)$ by elementary linear algebra. By performing long division, we obtain: $\Phi_5(x) = (x^2 +(1-t)x -t(1-t))p(x)+(t+t^2-t^3)x -t^2+t+1$, and since $t^2-t-1 =0$, $p(x)$ divides $x^5 -1=(x-1)\Phi_5(x)$. This indicates $X$ has order 5. \end{proof}
\emph{Remark}. Given an element $X$ in SL$(2,F)$ for any arbitrary field $F$, if $X$ has order $n$ then the possible orders for $\pi(X)$, where $\pi$: SL$(2,F)$ $\rightarrow$ PSL$(2,F)$ is the canonical projection defining PSL$(2,F)$ = SL$(2,F)$/ SL$(2,F)$ $\cap$ $\mathcal{Z}$, is either $n$ or $n/2$. ($\mathcal{Z}$ denotes the central scalar matrices over $F$). Thus, in this case, if $X$ has order 5, $\pi(X)$ is also of order 5. \

\begin{proof} (\emph{of theorem 2.4}) Consider the same matrices $M,N,W$ as in the proof of theorem 2.12, we define $R' = \begin{pmatrix} \alpha & \beta \\ \gamma &\delta \end{pmatrix}$. Since $4 |p-1$, the mutiplicative group $F_p^{*}$ contains an element of order 4 denoted as $i$. In addition, $p=\pm1$ mod $10$ implies that $5$ is a square mod $p$, thus, $\sqrt{5}$ makes sense in $F_p^*$. We denote $\mu_1 = \frac {1+\sqrt{5}}{2}$ and set $\gamma = - \mu_1$, $\beta=0$, $\alpha =i$, $\delta = -i$. This gives $R'= \begin{pmatrix} i&0\\-\mu_1 & -i\end{pmatrix}$ and $WR'= \begin{pmatrix} \mu_1&i\\i&0\end{pmatrix}$. By our construction, $R'$ and $WR'$ both have determinant 1 and $R'$ has order 4 since it is traceless. Also, since $\mu_1$ is a root of the equation $x^2 - x -1$, by proposition 2.5, $WR'$ is an element of order 5 in $SL(2,p)$. Thereby, $wr'$ has order 5 and $r'$ is an involution in $PSL(2,p)$ as noted earlier. Since $w$ and $r'$ have order 2, $r'w=w(wr')w,(wr')^{-1}=r'w, r'=(r'w)w$. This gives us the following representation:\\
 $\langle r',w\rangle=\langle w,wr'\rangle=\langle w,wr': w^2=(wr')^5=1, (wr')^{-1}=w(wr')w\rangle\cong D_{10}$ noncyclic, which indicates that $w,r'$ are distinct. 
 
\textbf{Claim 2.5.1}: $\langle wm,wr'\rangle \cong \langle wn,wr'\rangle \cong A_5$

\emph{Proof of claim}. $WR'WM=\begin{pmatrix} -x\mu_1+iy & y\mu_1 + ix \\ -ix & iy\end{pmatrix} \Rightarrow$ $Tr(WR'WM)=-x\mu_1 +2iy$. We set this trace value to be $\mu_1$. By Fricke's lemma, $Tr([WM,WR'])$ = $-4y^2-4ix\mu_1y-\mu_1^2y^2 -2$. Since $Tr(WM)=0$ and $Tr(WR')=\mu_1$, if we set this trace value to be 1, the main theorem in D. McCullough's paper provides $\pi(\langle WR', WM\rangle)\cong A_5$, where $\pi$ is the canonical projection from SL$(2,p)$ onto PSL$(2,p)$. Thus, our task is to show that there exists such $x,y$ satisfying the constraints. By using the relation $x^2 +y^2=-1$, we obtain the following quadratic equation: $-(1+\mu_1^{2})y^2 - (4ix\mu_1)y +3x^2=0 $. Solving for $y$ gives $y=-\frac{4i\mu_1\pm\sqrt{-4\mu_1^{2}+12}}{8}x$ = $\frac{-4i\mu_1\pm(\sqrt5-1)}{8}x$ (*). Inserting this relation to the equation $-x\mu_1 +2iy=\mu_1$ gives us $x = \pm(6+2\sqrt{5})$. As we notice beforehand, $\sqrt{5}$ makes sense in $F_p^{*}$. As the result, $\langle wm,wr' \rangle$'s congruence to $A_5$ is achieved. Likewise, if we fix our choice of one of the relations between $x$ and $y$ in (*) and pick the other for the relation between $z$ and $t$, the same procedure of finding distinct value for $z$ would also arrive at $\langle wn,wr'\rangle \cong A_5$. This proves our claim. 

In addition, since $r'$ is an involution in PSL$(2,p)$, by the same counting argument carried out in the proof of claim 2.3.2 and the fact that $4|A_5| < |G|$ for $p\geq 41$, the statements $\langle wm,wr' \rangle=\langle m,r',w\rangle$ and $\langle wn,wr'\rangle =\langle w,r',n\rangle$ are true (2). In fact, since we inherited the same construction for $W,M,N$, properness of $\langle wn,wm\rangle$ and $\langle n,m,w\rangle$ are imported into this proof. From this, the remaining task is to show that $(wn,wm,wr')$ is an irredundant generating sequence of $G$ and $w$ fails to replace any of the generators to give a new generating sequence. To see the irredundance of the sequence, we notice that $wn\notin \langle wm,wr'\rangle$; otherwise,by claim 2.4.1 and (2), $mn$ is an element in $\langle wm,wr'\rangle \cong A_5$. However, a straightforward calculation considering the trace of $MN$ would imply that $ord(\pi(MN)) \notin \{1,2,3,5\}$ for any possible pair of $(x,y)$ and $(z,t)$, which then leads to a contradiction. By symmetry, we also have $wm\notin\langle wn,wr'\rangle$. From this, if $wr'\in\langle wn,wm\rangle$, it follows that $\langle wm,wr'\rangle\cong A_5$ properly contained in $\langle wn,wm\rangle$, but this is impossible since $A_5$ is a maximal subgroup in PSL$(2,p)$ for $p=1$ mod $10$ and $\langle wm,wn\rangle$ is a proper subgroup. This result also implies $\langle wn,wm,wr'\rangle = G$ because it properly contains $\langle wm,wr'\rangle$. Last but not least, an entirely analogous justification done in the last paragraph for the proof of theorem 2.1 shows that $w$ cannot replace any of the generators to produce another generating sequence, which completes our proof.\end{proof}

\emph{Remark}. Note that, in theorem 2.4, for $G$ = PSL$(2,p^2)$ where $p=\pm$3 mod 10, since any quadratic equation splits in the field $F_{p^2}$, $\sqrt{5}$ is well-defined in $F_{p^2}$ as a root of $x^2 -5 =0$. Also, since $p^2 -1$ is always divisible by 4, the element $i$ of order 4 exists. As for the subgroup structure of $G$,  $A_5$ is maximal for $q=p^2$ where $p=\pm3$ mod $10$. From this, we now acquire all the necessary tools to create an entirely analogous construction by manipulating traces and orders of elements in SL$(2,p^2)$. In other words, by reproducing the same proof as in theorems 2.1 and 2.4, we reach the following result.\\\\ \textbf{Theorem 2.6}: For $G$=PSL$(2,p^2)$ with prime $p=\pm3$ mod $10$, $G$ fails RP if $m(G)=3$.\\\\ Thus far, we only devote our attention to the case when $m=3$ for PSL$(2,q)$ and when $m=4$, we know from Corollary 4.2 of \cite{nachman} that PSL$(2,p)$ satisfies RP. We may ask if there are other possible values for $m$ of PSL$(2,q)$ for arbitrary $q$. This question still remains unsolved. However, in Julius and Saxl' paper \cite{julius}, it is shown that the values 3 and 4 are exhaustive for $m$ of PSL$(2,p)$ where $p$ is any prime. Specifically, unless $p\equiv \pm$ 1 mod 8 or $p\equiv \pm$ 1 mod 10, $m$ is 3. Nachman later extended this result to $m=3$ if $p\not\equiv \pm$ 1 mod 10 except for PSL$(2,7)$ whose $m=4$ \cite{nachman}. In addition, Julius and Saxl established an upperbound for $m$ of PSL$(2,p^k)$, namely, $m\leq max(6,\pi(k) + 2)$ where $\pi(k)$ denotes the number of distinct prime divisors of $k$. Further classification for $m$ of PSL$(2,q)$ regarding the so-called subgroups in \emph{general position}, which will be introduced in section 3, is also provided in Theorem 7 of \cite{julius}.\\

\emph{Remark}. It is important to notice that since PSL$(2,q)$ are (2,3)-generated for any $q\neq9$, that is the group can be generated by an involution and an element of order 3, it can also be generated by 3 conjugate involutions \cite{mac69}. This sets a sharp lower bound for $m$ of PSL$(2,q)$ where $q\neq9$. A stronger result, which initially sets out to address the question whether PSL$(n,q)$ can be generated by three involutions, two of which commute, also arrives at the same lower bound for $m$ of PSL$(2,q)$ (see \cite{cherkassoff}). This raises another question about how often the value of $m$ equals to 3. Again, for any arbitrary $q$, this still remains an open problem. However, for the case where $q$ is prime, an interesting result in Jambor's paper \cite{jambor} states the following.\\ \newpage
\textbf{Theorem 2.7} (S.Jambor): PSL$(2,p)$ contains an irredundant generating sequence of length 4 if and only if $p\in\{7,11,19,31\}$. More precisely, up to automorphism, there are 2 irredundant generating sequences of length 4 for PSL$(2,7)$, 14 for PSL$(2,11)$, 3 for PSL$(2,11)$ and 1 for PSL$(2,31)$.\\

\emph{Remark}. This classification along with Whiston and Saxl's result implies that $m$ = 3 for all primes except for the listed. Yet, whether the ubiquity of the value $3$ for $m$ holds for non-prime $q$ is subject to further investigation.\\
\begin{indent}\end{indent}Considering higher values of $m$, say 4, one may hope for the satisfaction of the replacement property for PSL$(2,q)$, possibly with $q$ a higher prime power since this is certainly the case for PSL$(2,p)$. Up to this point, we only utilize the definition of the replacement property to form explicit sequences that fail. To show the property holds, we need more information connecting generating sequences with appropriate sequences of maximal subgroups. We will employ such relation to obtain the desired result for a class of irredundant generating sequences, namely, ones made of only involutions.      
%--------------------------------------------

\section{Involutions and RP-satisfaction in $\mathrm{PSL}(2,p^2)$}
\textbf{Definition}: A collection of subgroups $\{H_i\leq G\}_{i\in I}$ with index set $I$ is said to be in \emph{general position} if for every $j\in I$, the intersection $\bigcap_{i\in I} H_i$ is properly contained in $\bigcap_{i\in I - \{j\}}H_i$.

\vspace{0.1cm}This definition arises naturally from an analogous example from linear algebra: consider an $n$-dimensional vector space $V$ with basis $\mathcal{B}$=$\{v_i: 1\leq i\leq n\}$, subspaces $W_i$=$span(v_j:j\neq i)$ are in general positions.  

To relate the definition with the study of finite generating sets of group, suppose we are given an irredundant sequence $s=(g_1,...,g_k)$ of a group $G$, for each $i$, let $H_i$ be the subgroup generated by the generators $g_j$ where $j\neq i$, i.e $H_i = \langle g_j:i\neq j\rangle$. At this point, it readily follows that the subgroups $H_i$'s are in general position since their intersection does not contain any generator $g_i$, whereas intersection of any $k-1$ members of the collection must contain the subgroup generated by exactly one generator. Yet, the same argument can be made to obtain a collection of maximal subgroups in general positions. Since each $H_i$ is a proper subgroup of $G$, there exists a maximal subgroup $M_i$ strictly containing $H_i$, therefore, containing generators $g_j$ for $j\neq i$. This gives rise to a dimension-like invariant for finite groups, whose properties and further theoretical/computational application are introduced and motivated in \cite{Ravi}. 
   
For the interest of this section, we note several important relations between maximal subgroups in general position and RP-satisfaction.\\

\begin{noindent}\end{noindent}\textbf{Proposition 3.1}(D. Collins and R.K.
Dennis)\ : Let $s=(g_1,...,g_k)$ be an irredundant generating sequence for
a group $G$. If every sequence of corresponding maximal subgroups in general
positions, say $(M_1,...,M_k)$, intersects trivially, $s$ satisfies the replacement
property.
\begin{proof}
Suppose $s$ fails RP, that is, there exists a nontrivial element g in $G$ so that for each slot $i$-th in sequence $(g_1,...,g_k)$, $g$ cannot replace $g_i$ to give a new generating sequence, i.e $H_i=\langle g_1,...,g,...,g_k\rangle \lneq G$. Thus, there exists a family of maximal subgroups $M_i$ properly containing $H_i$. By definition, $(M_i)^k_{i=1}$ is a corresponding sequence of maximal subgroups in general position and $\bigcap^k_{i=1} M_i$ contains $g$, hence a contradiction.     
\end{proof}

\begin{noindent}\end{noindent}\textbf{Proposition 3.2}: Let $G$ be a finite
group, for any $k\leq m(G)$, let $s=(g_1,...,g_k)$ an irredundant generating
sequence for $G$. For any corresponding collection of maximal subgroups in
general positions $(M_1,...,M_k)$, if there exists $r$ in $\{1,...,k\}$ such
that the following hold:\begin{enumerate}

\item $M_r=\langle g_i:i\neq r\rangle$,
\item $m(M_r)=k-1$,
\item $M_r$ satisfies the replacement property,
\end{enumerate}
then $s$ satisfies the replacement property.

\begin{proof} 
Given such a collection of maximal subgroups, without loss of generality,
we let $r=k$. Let $M_i'=M_i\cap M_k$, where $ i:1 \leq i \leq k-1$, since
the $M_i$ 's are in general position, $M'_i\geq\langle g_j\rangle$ iff $j\neq
i, k$ . This implies  the $M'_i$ 's are in general position with respect
to $(g_1,...,g_{k-1})$. Also, since $M'_i  \lneq M_k$, there exists $N_i$
maximal subgroups of $M_k$ containing $M_i'$. Thus, $N_i$ 's are maximal
subgroups of $M_k$ in general position. Now, we have $\bigcap_{i=1}^kM_i$= $\bigcap_{i=1}^{k-1}M_i'\leq\bigcap_{i=1}^nN_i$=$\{ e\}$, where the last equality
comes from the fact that $M_k$ satisfies the replacement property. By proposition
3.1, we have $s$ satisfies the replacement property.
\end{proof}

This relation establishes an essential criterion on the maximal subgroups, which allows the satisfaction of RP of irredundant generating sequence of order-2 elements.\\     

\begin{noindent}\end{noindent}\textbf{Theorem 3.3}: For $p$ odd prime,
$ p > 5$, $G = \mathrm{PSL}(2,p^2)$, $m(G)\geq
4 $, any length-4 irredundant generating sequence of involutions satisfies RP.\\

The existence of more than two involutions in an irredundant generating sequence for $G$ provides an implicit insight into the subgroup structure of $G$, namely the dihedral subgroups generated by the generators with appropriate orders imply the uniqueness of subgroups containing them. We shall see why this is true via the following propositions.\\ 

\begin{noindent} \textbf{Proposition 3.4}: Let $D_{2p}$ be a dihedral group of order $2p$ in $G$=PGL$(2,p^2)$, for any odd prime $p$, $N_G(D_{2p})$ is isomorphic to $C_p \rtimes C_{p-1}$.\end{noindent}
\begin{proof}
 We observe that there is a unique element on the projective line in the natural actions
of PGL$(2,p^2)$ of degree $p^2+1$, that is fixed by $D_{2p}$. This means that the element is also fixed by $N_G(D_{2p})$ since  the normalizer of a subgroup contains the subgroup itself. Thus, $N_G(D_{2p})$ is contained in the point stabilizer, which, by a standard result in group theory, is an affine group of form $M\rtimes N$, where $M\rtimes N = F_{p^2}^{+} \rtimes_\rho F_{p^2}^{*}$. The map $\rho$ is the action of $N$ on $M$ via field multiplication. This gives us the group operation in $N_G(D_{2p})$. Now, we regard $D_{2p}$ as $D:=\{(a,b):b=\pm1, a\in F_{p}\}$. The task is to deduce $N_G(D_{2p}) \cong F_p^{+} \rtimes F_p^{*}$. To do this, it suffices to show that for any element $(m,n)$ in the normalizer of $D$, $n$ is an element in the multiplicative group of $p-1$ elements and $m$ is in the additive group of $p$ elements. To see why this is true, let us consider $(m,n)$'s interaction with the generators of $D$, namely, $(0,-1),(1,1)$, which gives: $(m,n)(0,-1)=(a,b)(m,n)$ for some $(a,b)\in D$ \
$\Rightarrow (m,-n)=(a+bm,bn)$\
$\Rightarrow b=-1$ and $m=2^{-1}a$, which is in $F_p$ where both factors make sense for $p$ sufficiently large.
Also, $(0,-1)(m,n)=(m,n)(a',b')$ for some $(a',b')\in D$
$\Rightarrow (-m,-n)=(m+na',nb')$ \ 
$\Rightarrow b=-1, -2m=na'$. This means that if $n=-2m(a')^{-1}$, where $m,a$ cannot be zero, which implies $n\in F_p^{*}$. For the case, $m$ and $a'$ are zero, we have $(0,n)(1,1)=(a'',b'')(0,n)$ for some $(a'',b'')$ in $D$ \
$\Rightarrow (n,n)=(a'', nb'')$\ 
$\Rightarrow n=a''$. Since $n \in F_{p^2}^{*}$ and $a''\in F_p$, $n \in F_p^{*}$.\end{proof}\vspace{0.10cm}

\begin{noindent}\end{noindent}\textbf{Proposition 3.5}: Let $D_{2p}$ be a dihedral subgroup of order $2p$ of $G$ = PSL$(2,p^2)$,
where $p$= $\pm$ 1 mod 4, there exists a unique subgroup $H$ of $G$ such
that $H$ $\cong$ PSL$(2,p)$, and $H$ contains $D_{2p}$.
\begin{proof}By Dickson's classification of subgroups of PSL$(2,p)$, we first
observe that since the order of PSL$(2,p)$ is either   $p(p-1)(p+1)$ or \
    $\frac{1}{2}p(p-1)(p+1)$, for PSL$(2,p)$ to contain a dihedral group
of order $2p$, $p-1$\ or $p+1$ are divisible by 4, hence the condition on
$p$. Also, such existence of $D_{2p}$ is also guaranteed due to the divisibility
of   $p^2-1$  by 8, hence by 4. The number of subgroups  isomorphic to PSL$(2,p)$
in PSL$(2,p^2)$ is given by the index of PGL$(2,p)$ in PGL$(2,p^2)$, namely,
$\frac{|\mathrm{PGL}(2,p^2)|}{|\mathrm{PGL}(2,p)|}$      $ (*)$.

Since all the dihedral groups of order $2p$\ are conjugate of one another
in PGL$(2,p^2)$, by proposition 3.4, we obtain $\frac{|\mathrm{PGL}(2,p^2)|}{(p-1)p}$
subgroups isomorphic to $D_{2p}$. Likewise, the number of  subgroups isomorphic
to $D_{2p}$  in $H$, where $H$ is any subgroup of $G$ isomorphic to PSL$(2,p)$
is given by $\frac{|\mathrm{PGL}(2,p)|}{(p-1)p}$ $(**)$. We notice that the product
of $(*)$ and $(**)$ is exactly the number of dihedral subgroups of order $2p$
in PSL$(2,p^2)$, which implies that the any each subgroup isomorphic to $D_{2p}$
uniquely determines a subgroup in PSL$(2,p^2)$, that is isomorphic to
PSL$(2,p)$. This is sufficient to complete the proof. \end{proof} 
By using the same strategy and Dickson's classification of subgroup structure of PSL$(2,p^k)$, we obtain the following similar result about the unique determination of subgroups in PSL$(2,q)$ containing a dihedral substructure.\\\\\textbf{Proposition 3.6}: If $D_{2n}$ is a dihedral subgroup of order $2n$ of $G$=PSL$(2,p^2)$ where $2n|p\pm1$, there exists a unique subgroup $H$ of $G$ isomorphic to PSL$(2,p)$ containing $D_{2n}$. 
\begin{proof}
By \cite{dickson}, there are $\frac{p(p^2-1)}{2n(2,p-1)}$=$\frac {|\mathrm{PGL}(2,p)|}{4n}$ dihedral subgroup $D_{2n}$ in PSL$(2,p)$. Likewise, there are $\frac{|\mathrm{PGL}(2,p^2)|}{4n}$ copies of $D_{2n}$ in PSL$(2,p^2)$ since $2n|p^2-1$. By conducting the same counting argument as in the previous proof, we obtain the desired result. \end{proof}\vspace{0.10cm}
\begin{noindent}\end{noindent}\textbf{Proposition 3.7}: Suppose $H_1$,
$H_2$ are subfield subgroups of PSL$(2,p^2)$, if the intersection of $H_1
$ and $H_2$ contains a dihedral subgroup of order $2p$, where $p\pm 1$ is
divisible by 4, $H_1\cap H_2$ is also a subfield subgroup. In addition, the
same result holds if the intersection contains a $D_{2n}$, where $2n$ divides
$p\pm1$.  
\begin{proof}
Given such $H_1$, $H_2$ and $D_{2p}$ dihedral  subgroup of $H_1 \cap H_2$,
$H_i$ 's are covering groups of PSL$(2,p)$, where PSL$(2,p)$ can be  embedded
into $H_i$ via isomorphism with a subgroup of $H_i$ ($i=1,2$). From this,
we have the following cases: 

\textbf{Case 3.7.1}: If $H_1$ $\cong$ PSL$(2,p)$, $H_2$ $\cong$ PSL$(2,p)$, we claim
that $H_1\cap H_2$ $\cong$ PSL$(2,p)$.\

\emph{Proof}. Since $H_1, H_2$\ both contains $D_{2p}$, by the uniqueness
part in proposition 3.5, we must have $H_1=H_2$, hence the claim.

\textbf{Case 3.7.2}: If $H_1$ $\cong$ PSL$(2,p)$, $H_2$ $\cong$ PGL$(2,p)$, we claim
that $H_1 \cap H_2$ $\cong$ PSL$(2,p)$. \
  
\textit{Proof}. There exists $K_1$ unique subgroup
of $H_2$ isomorphic to PSL$(2,p)$ containing $D_{2p}$. By the uniqueness
part of proposition 3.5 again, we have $K_1 = H_1$, which implies $H_1\cap
H_2=H_1$ $\cong$ PSL$(2,p)$.

\textbf{Case 3.7.3}: If $H_1$ $\cong$ PGL$(2,p)$, $H_2$ $\cong$ PGL$(2,p)$, we claim
that $H_1 \cap H_2$ $\cong$ PGL$(2,p)$ \

\emph{Proof}. By the same reasoning,  there exist $K_1,K_2,$ respectively,
unique subgroups of $H_1, H_2$ isomorphic to PSL$(2,p)$ containing $D_{2p}$
, and $K_1=K_2$. Also, we note that $K_1=K_{2}$ is contained in a unique
subgroup of PSL$(2,p^2)$ isomorphic to PGL$(2,p)$. This indicates that $H_1=H_2$,
hence the result.

With proposition 3.6, a similar proof applies for the case where the intersection contains a dihedral subgroup of order $2n$ satisfying the hypothesis.  \end{proof}
\begin{proof}\emph{(of theorem 3.3) }Consider such a generating sequence
$s=(g_1,g_2,g_3,g_4)$, where the $g_i$'s are involutions, we let $M_1,M_2,M_3,M_4$
be an associated maximal subgroups in general positions, where $M_i\geq\langle
g_j:j\neq i \rangle $ and $i=1,2,3,4$. By Dickson's classification of maximal
subgroups for PSL$(2,q)$, these maximal subgroups are of the classes: \\
$(1)$ Stabiliser of a point in $P_1(q)$  \\
$(2)$ $D_{q-1}$ for $q$ odd and $D_{2(q-1)}$ for $q$ even\\
$(3)$ $D_{q+1}$ for $q$ odd and $D_{2(q+1)}$ for $q$ even\\
$(4)$ subfield subgroups PSL$(2,q_1).a$ for $a \leq 2$\\
$(5)$ $A_4.a$, where $q=p\geq5$ , $a\leq 2$\\
$(6)$ $A_5$ where $q=\pm1$ mod $10  $. 

We notice that there no maximal subgroup of type (5) since $q$ is non-prime.
$q_1$ is defined to be $p^{\frac{r}{p_1}}$, where $r$ is the power of $p$,
  $p_1$ is any prime divisor of $r$. In our case, $r=2$, hence $q_1=p $ .
Thus, the only subfield subgroups of class (4) are isomorphic to PSL$(2,p)$
or PGL$(2,p)$. Also, by proposition 2 in Whiston and Saxl's paper, we can
only have at most 3 $M_i$ 's are of the first three classes. In fact, there
can only be at most two $M_i$'s are of the first three classes since $m(G)\geq
4$. This shows that there at least two $M_i$'s that are of class (4) or (6).

\textbf{Claim 3.3.1}: There can only be at most one $M_i$, that is a subfield subgroup.

\textit{Proof of claim}. Suppose we have at least two $M_i$ 's that are subfield
subgroups. Without loss of generality, let $M_1\gneq \langle g_2,g_3,g_4\rangle
$, $M_2\gneq \langle g_1,g_3,g_4\rangle$, then $M_1\cap M_2 \geq \langle g_3,g_4
\rangle $. Since $g_3,g_4$ are involutions, $M_1\cap M_2$\ contains a dihedral
subgroup of order $2n$ for some $n$, but we also note earlier that the $M_i$'s
can only be $PGL(2,p)$ or $PSL(2,p)$ . As noted in previous propositions,
to contain such a dihedral subgroup, we either have $2n=2p$ where $p\equiv$ 1 mod 4 or $2n$ divides $p\pm 1$.
By proposition 3.7, the intersection of $M_1$ and $M_2$ is a subfield subgroup.
Also, by the proof of proposition 3.7 itself, there are only three possibilities
for this intersection. However, since $M_1\cap M_2 \neq M_1$ or $M_2$,
by checking each case, we always have $M_1\cap M_2 \cong M_1$ or $M_2$,
which is impossible since the cardinality of the intersection is strictly
less than the cardinalities of the $M_i$'s. This shows the claim is true.\

As the result, there must be at least one of the maximal subgroups in general
position that is  $A_5$, say $M_3$ . By lemma 3.2 in B. Nachman, any length-3
irredundant sequence in $A_5$\ must generate $A_5$. This means that since
$M_3$ contains a subgroup generated by an irredundant sequence of 3 involutions,
such sequence generates $M_3$. Also, since $A_5$ satisfies the replacement
property and $m(A_5)=3$, by proposition 3.2, $s$\ satisfies the replacement
property. This completes our proof. 
\end{proof}
From the proof of claim 3.3.1, it is an entirely identical situation when the length of $s$ is larger than 4. That is, from the result of Whiston and Saxl, one can apply the same argument to give the following statement.\\

\begin{noindent}\end{noindent}\textbf{Corollary 3.8}: Suppose prime $ p > 5$, $G$= PSL$(2,p^2)$ and 
$m(G)\geq4$, given an irredundant generating sequence of involutions of length at least 4, for any associated collection of maximal subgroups in general position, there is at least one member of such collection isomorphic to $A_5$.
\begin{proof}
 We denote $m=m(G)>4$. The case where the length of sequence equal to 4 is treated by the previous proof. Thus, for an irredundant generating sequence of order-2 elements $s=(g_1,...,g_k)$, where $4< k\leq m$, let us have  $\{M_i\}_{i=1}^k$ be an associated collection of maximal subgroups in general positions. Again, since $m>3$, there can only be at most two members of the collection belonging to the first three classes. This setting and the proof of claim 3.3.1 give us at least 2 maximal subgroups in general position as copies of $A_5$.\end{proof}

\emph{Remark}. This leads to another interesting result. Let us consider the same values $m$, $k$ and sequence $s$ as in the previous proof. The corollary asserts that $A_5$ is always a member of any corresponding maximal subgroup in general position. Without loss of generality, we set $M_k$ as $A_5$. By definition, we have $M_k > \langle g_1,g_2, ...,g_{k-1}\rangle$. Since the $g_i$'s are irredundant and any length-3 irredundant sequence in $A_5$ generates $A_5$, $(g_1, g_2,...,g_{k-1})$ irredundantly generates $A_5$. This is a contradiction because $m(A_5)=3$ while $k-1\geq 4$. Thus, we arrive at the following corollary.\\

\begin{noindent}\end{noindent}\textbf{Corollary 3.9}: For prime $ p > 5$, $G$= PSL$(2,p^2)$ and 
$m(G)\geq4$, any irredundant generating sequence of involutions has length at most 4.  

\emph{Remark}. As noted earlier, this kind of result is expected; however, it is not straightforward to be proven without considering the intersection of possible maximal subgroups that are subfields in $G$. In fact, we are relying on the existence of dihedral substructures within the corresponding maximal subgroups in general position and their intersections to eliminate certain classes of subgroup description and to guarantee existence of subgroup such as $A_5$. 

\section{Acknowledgment}
The author wishes to express his utmost appreciation towards Professor R.K Dennis and Ravi Fernando, who provided insightful discussions and guidance upon this topic of finite group theory, which allow the possibility of this paper. In addition, the author is particularly grateful for Benjamin Nachman's paper, which provide the strategies for theorems 2.1 and 2.4.

\end{document}